\documentclass[nointlimits,11pt,oneside]{amsart}
\usepackage{amssymb,cases,enumerate}
\usepackage{color}
\usepackage{hyperref,mathscinet}
\usepackage{natbib}
\usepackage{mathscinet}
\usepackage{hyphenat}

\usepackage[scr=boondoxo,scrscaled=1.05]{mathalfa}

\usepackage[%
	a4paper,
	total={16cm,23cm},
	top=3cm,
	marginparsep=2pt]
{geometry}

\newcommand{\N}{\mathbb{N}}
\newcommand{\R}{\mathbb{R}}

\newcommand{\M}{\mathcal M}
\newcommand{\Mpl}{\M_+}
\renewcommand{\d}{\,d}
\renewcommand{\b}{{\mathscr{b}}}

\DeclareMathOperator*{\esssup}{ess\,sup}
\DeclareMathOperator{\sgn}{sgn}

\theoremstyle{plain}
\newtheorem{theorem}{Theorem}[section]

\theoremstyle{definition}
\newtheorem{remark}[theorem]{Remark}
\newtheorem{remarks}[theorem]{Remarks}
\newtheorem{definition}[theorem]{Definition}

\newtheorem*{conventions}{Conventions}

\numberwithin{equation}{section}

\usepackage[textsize=tiny,textwidth=2cm,colorinlistoftodos]{todonotes}

\makeatletter
\def\paragraph{\bigskip\@startsection{paragraph}{4}%
  \z@\z@{-\fontdimen2\font}%
  {\normalfont\bfseries}}
\makeatother

\title{Reduction principle for Gaussian $K$-inequality}
\author{Sergi Baena-Miret, Amiran Gogatishvili, Zden\v ek Mihula and Lubo\v s Pick}

\address{Sergi Baena-Miret,
Department of Mathematics and Computer Science, University of Barcelona, Gran Via de les Corts Catalanes 585, Barcelona, 08007, Spain}
\email{sergibaena@ub.edu}
\urladdr{0000-0002-8423-5098}

\address{Amiran Gogatishvili,
 Institute of Mathematics of the
 Czech Academy of Sciences,
 \v Zitn\'a~25,
 115~67 Praha~1,
 Czech Republic \newline
 L. N. Gumilyov Eurasian National University,
5 Munaytpasov St., 010008 Nur-Sultan, Kazakhstan}
\email{gogatish@math.cas.cz}
\urladdr{0000-0003-3459-0355}

\address{Zden\v ek Mihula, Czech Technical University in Prague, Faculty of Electrical Engineering, Department of Mathematics, Technick\'a~2, 166~27 Praha~6, Czech Republic -- AND -- Department of Mathematical Analysis, Faculty of Mathematics and Physics, Charles University, Sokolovsk\'a~83, 186~75 Praha~8,	Czech Republic}
\email{mihulzde@fel.cvut.cz \& mihulaz@karlin.mff.cuni.cz}
\urladdr{0000-0001-6962-7635}

\address{Lubo\v s Pick,
 Department of Mathematical Analysis,
	Faculty of Mathematics and Physics,
	Charles University,
	Sokolovsk\'a~83,
	186~75 Praha~8,
	Czech Republic}
\email{pick@karlin.mff.cuni.cz}
\urladdr{0000-0002-3584-1454}

\begin{document}

\date{\today}



\thanks{This research was supported in part by the grant P201-18-00580S of the Czech Science Foundation and by the Danube Region Grant no.~8X2043. The research of  A.~Gogatishvili was also supported by Czech Academy of Sciences RVO: 67985840. The research of Z.~Mihula was supported by the project OPVVV CAAS CZ.02.1.01/0.0/0.0/16\_019/0000778 and by the grant SVV-2020-260583. The research of S.~Baena-Miret was supported by the grant 2017SGR358.}

\begin{abstract}
We study interpolation properties of operators (not necessarily linear) which satisfy a specific $K$-inequality corresponding to endpoints defined in terms of Orlicz--Karamata spaces modeled upon the example of the Gaussian--Sobolev embedding. We prove a reduction principle for a fairly wide class of such operators.
\end{abstract}

\maketitle

\setcitestyle{numbers}
\bibliographystyle{abbrv}

\section{Introduction}

The principal motivation for our research is investigate the applicability of interpolation techniques, in particular the $K$-method, to sharp Gaussian--Sobolev embeddings, or, more generally, Boltzman--Sobolev embeddings. Such an approach was successfully applied earlier for example to Euclidean--Sobolev embeddings (\cite{KePi:06}), boundary
trace embeddings (\cite{CiKePi:08}), or to a wide variety of classical operators of harmonic analysis (\cite{EdMiMuPi:20}). The method can be outlined as follows: we begin with two sharp endpoint estimates from which an inequality between corresponding $K$-functionals is derived (we will refer to this step as a
\textit{$K$-inequality}). The $K$-inequality typically gives a pointwise comparison of certain operators involving nonincreasing rearrangements of images of an operator to those involving the functions themselves. This inequality is then corroborated using ideas involving some modification of the Hardy--Littlewood--P\'olya relation. In case of Sobolev embeddings, some result in the spirit of the DeVore--Scherer theorem is used at the end.

In~\cite{KePi:06},~\cite{CiKePi:08} and~\cite{EdMiMuPi:20}, this approach worked very well. On the other hand, for example in~\cite{CP:09}, where sharp Gaussian--Sobolev embeddings were established, interpolation methods were not used. Instead, the optimal embeddings were derived from an appropriate isoperimetric inequality. This step required a symmetrization argument exploiting a general P\'olya--Szeg\H{o} principle on the decrease of rearrangement-invariant norms of the gradient of Sobolev functions in the Gauss space, extending earlier results of \cite{Eh:84} and \cite{Ta:92}. The proof relied upon the
Gaussian isoperimetric inequality by Borell \cite{Bo:75} which gives an explicit description of the isopertimetric profile of $\mathbb R^n$ endowed with the probability Gaussian measure. A serious disadvantage of this technique consists in the fact that it works only for the first-order embeddings. This shortcoming was later overcome by establishing higher-order results using sharp iteration methods (\cite{CiPiSl:15}). Thanks to those results, sharp function spaces appearing in such embeddings are known, at least in the rearrangement-invariant environment.

In the light of the described situation, it would clearly be of interest to investigate the very existence of an operator (or operators) whose boundedness between a given pair of rearrangement-invariant function spaces would guarantee that every operator satisfying the $K$-inequality corresponding to specific pairs of endpoint spaces, modeled upon the example of endpoint spaces appropriate for Gaussian--Sobolev embeddings, is bounded between that pair. This idea is to some extent connected with the classical result of
Calder\'on (\cite{Ca:66}). However, in this paper, we are not so much interested in characterizing Calder\'on couples, but instead we aim at nailing down those pairs of function spaces for which the corresponding $K$-inequality always guarantees the boundedness of operators. Although our research was originally motivated by Gaussian--Sobolev embeddings, operators having similar endpoint behavior appear also in other circumstances, for instance in studying problems appearing in Gaussian harmonic analysis (\cite{Ur:19}).

A prototypical example, motivated by the Gaussian--Sobolev embeddings, of endpoint behavior that we have in our mind is that of an operator $T$ satisfying the $K$-inequality
\begin{equation}\label{E:endpoints-of-S}
K(Tf, t; L\sqrt{\log L}, e^{L^2}) \lesssim K(f,t;L^1,L^\infty) \quad \text{for every $t\in(0,1)$}
\end{equation}
with a multiplicative constant independent of $f$. It will be useful to notice that both of the spaces on the left-hand side are the classical Orlicz spaces of either logarithmic or exponential type, sometimes also called \textit{Zygmund classes}. It it also important to recall that these spaces are neither Lebesgue spaces nor two-parameter Lorentz spaces, which makes their study through interpolation techniques considerably difficult. On the other hand, they are special cases of the Lorentz--Zygmund spaces (\cite{BeRu:80}), and also of the yet more general Lorentz--Karamata spaces, based on the so-called \textit{slowly-varying functions}. These spaces were first introduced in~\cite{EdKePi:00} and then treated by many authors (see e.g.~\cite{Bat:18,GoNeOp:04,GoOpTr:05,Pe:21}). It might be useful to note that, in the notation of Lorentz--Zygmund spaces~\cite{BeRu:80},~\eqref{E:endpoints-of-S} reads as
\begin{equation}\label{E:endpoints-of-S-lz}
K(Tf, t; L^{1,1;\frac{1}{2}}, L^{\infty,\infty;-\frac{1}{2}}) \lesssim K(f,t;L^1,L^\infty) \quad \text{for every $t\in(0,1)$}.
\end{equation}
It turns out that the principal property of every operator $T$ satisfying \eqref{E:endpoints-of-S-lz} is the validity of
\begin{equation}\label{E:K-inequality-gauss}
    \int_{0}^{t}\frac{(Tf)^*(s)}{\sqrt{\log \frac{e}{s}}}\,ds \lesssim \int_{0}^{t}\frac{f^*(s\log \frac{e}{\sqrt{s}})}{\log \frac{e}{s}}\,ds \quad \text{for every $t\in(0,1)$}.
\end{equation}
We shall use a far more general form of this inequality as a point of departure. Namely, for $p\in(0, \infty)$ and for a pair $(\b_1,\b_2)$ of slowly varying functions, we will consider operators $T$ satisfying
\begin{equation*}
\int_0^t \left[(Tf)^*(s)\b_1(s)\right]^p \, ds \lesssim \int_0^t \left[f^*\big(\sigma^{-1}(s^{1/p})^p\big) \b_1(s) \b_2(s)^{-1}\right]^p \, ds
\end{equation*}
\noindent for every $f\in L^p(R,\mu)$ and $t\in(0,1)$, where
\begin{equation*}
\sigma\colon[0,1]\to[0,1]
\end{equation*}
is the increasing, bijective function such that
\begin{equation*}
    t^p = \frac 1C\int_0^{\sigma(t)^p}\left[\b_1(s) \b_2(s)^{-1}\right]^p\d{s}\quad\text{for every $t\in[0,1]$}
\end{equation*}
for an appropriate constant $C$. Motivated by the principal inspiration and motivation, we shall call such operators \emph{$(p,\b_1,\b_2)$-gaussible}.

Our main result is Theorem~\ref{gaussible:calderon_thm:restated}, complemented with Theorem~\ref{gaussible:Usigma_bdd_implies_S_bdd}, below. It gives several characterizations of boundedness of every $(p,\b_1,\b_2)$-gaussible operator from $X$ to $Y$, where $X,Y$ is a prescribed pair of rearrangement-invariant spaces over two (possibly different) nonatomic measure spaces of measure 1. Of course, our choice of the value 1 is made only for technical convenience, and is immaterial as simple modifications can be used to extend the results to any nonatomic finite measure space. For the particular case corresponding to \eqref{E:endpoints-of-S}, Theorems~\ref{gaussible:calderon_thm:restated} and~\ref{gaussible:Usigma_bdd_implies_S_bdd} yield (among other results) that for a given pair of rearrangement-invariant spaces $X$ and $Y$, the following three statements are equivalent:
\begin{enumerate}[(i)]
    \item every operator $T$ satisfying~\eqref{E:K-inequality-gauss} is bounded from $X$ to $Y$,

    \item every operator $T$ satisfying \eqref{E:endpoints-of-S} is bounded from $X$ to $Y$,

    \item the operator $U$ defined by $Uf(s)=f^{*}(s\log\frac e{\sqrt{s}})\sqrt{\log\frac{e}{s}}$ for $s\in(0,1)$ and every suitable $f$ is bounded from $X(0,1)$ to $Y(0,1)$, where $X(0,1)$ and $Y(0,1)$ are the representation spaces of $X$ and $Y$ in the classical Lorentz--Luxemburg sense.
\end{enumerate}
It is worth noticing that the operator $U$ in (iii) is far away from being quasilinear, let alone linear.

To provide the interested reader with some useful information, we shall now describe the motivation and what lies at the root of \eqref{E:endpoints-of-S-lz} in more detail. The story begins with the seminal paper~\cite{Gr:75} of L. Gross, who established the first of Gaussian--Sobolev embeddings and also pointed out its importance. In the study of quantum fields and hypercontractivity semigroups, one often needs semigroup estimates, which can be equivalently described in terms of inequalities of Sobolev type in infinitely many variables (see, for instance,~\cite{Ne:73} and the references therein). In~\cite{Gr:75}, the major problem occurring in attempts to generalize classical Sobolev embeddings to cases of infinitely many variables (recall that the Lebesgue measure does not make sense for infinitely many variables) was solved by replacing the Lebesgue measure by the Gaussian probability measure in $\mathbb R^n$, $n\ge 1$, having the density
\begin{equation*}
    d\gamma_n(x) = (2\pi)^{-\frac{n}{2}}e^{\frac{-|x|^2}{2}}\,dx,
\end{equation*}
and requiring the embedding constants to be independent of the dimension. Since $\gamma_n(\mathbb R^n)=1$ for every $n\in\N$, taking limit as $n\to\infty$ makes sense. It should be mentioned though that another very important question was settled in the same paper, a question concerned with the comparison of integrability of the gradient of a scalar function of several variables with the integrability of the function itself. While in the Euclidean environment there always is a huge gain in integrability, expressible by change of certain power, no such thing is available in the Gaussian setting. Typically, if $\nabla u\in L^{p}(\mathbb R^{n},dx)$ for some $p\in[1,n)$, then $u\in L^{\frac{np}{n-p}}(\mathbb R^{n},dx)$, in which $dx$ stands for the $n$-dimensional Lebesgue measure, and, of course, $\frac{np}{n-p}>p$. But with $n\to\infty$ one has $\frac{np}{n-p}\to p$, so there is a good chance that the gain will be lost. However, L. Gross discovered that there still is some gain, albeit only of a logarithmic, rather than power, nature. Namely, he proved that if a function $u$ satisfies $\nabla u\in L^{2}(\mathbb R^n,\gamma_n)$ and is suitably normalized (for example when its median is zero), then $u$ itself belongs to a slightly ``better'' space (better means smaller here), namely $L^{2}\log L(\mathbb R^n,\gamma_n)$. A more precise formulation of this inequality reads as follows:
\begin{equation*}
    \|u-u_{\gamma_n}\|_{L^{2}\log L(\mathbb R^n,\gamma_n)} \le C \|\nabla u\|_{L^2(\mathbb R^n,\gamma_n)},
\end{equation*}
where $C$ depends on neither $u$ nor $n$, and $u_{\gamma_n}$ denotes the integral mean of $u$, that is
\begin{equation*}
    u_{\gamma_n} = \int_{\mathbb R^{n}} u(x)\,d\gamma_n(x).
\end{equation*}

The discovery of L. Gross paved the way for extensive research of all kinds. His result has been extended, generalized and modified many times, and simple proofs occurred as well as new applications. In~\cite{CP:09}, a comprehensive treatment of sharp Gaussian--Sobolev embeddings of the form
\begin{equation}\label{E:gaussian-sobolev-ri-embedding}
    \|u-u_{\gamma_n}\|_{Y(\mathbb R^n,\gamma_n)} \le C \|\nabla u\|_{X(\mathbb R^n,\gamma_n)}
\end{equation}
was carried out, in which $X$ and $Y$ are general rearrangement-invariant spaces. The focus has been on the ``optimality'' of the function spaces involved. One of the most important discoveries of~\cite{CP:09} was that the Gaussian--Sobolev embedding can be equivalently described by the action of an operator acting on functions of a single variable, providing thus a considerable simplification of the problem in hand. Namely, it was shown in~\cite[Theorem 3.1]{CP:09} that the inequality~\eqref{E:gaussian-sobolev-ri-embedding} is equivalent to the boundedness of the operator $S$ defined as
\begin{equation*}
    Sg(t) = \int_{t}^{1}\frac{g(s)}{s\sqrt{\log\frac{e}{s}}}\,ds
\end{equation*}
for suitable functions $g\colon (0,1)\to\mathbb R$ and every $t\in(0,1)$ from $X(0,1)$ to $Y(0,1)$. This result is usually called a \textit{reduction principle}. The operator $S$ is known to satisfy
\begin{equation}\label{E:endpoints-of-S-restricted}
    \begin{split}
        &S\colon L^1(0,1)\to L(\log L)^{\frac{1}{2}}(0,1),
            \\
        &S\colon L^\infty(0,1)\to \exp L^2(0,1),
    \end{split}
\end{equation}
and, interestingly, this ``endpoint behavior'' is shared also by the operator $U$. We shall, however, prove as a particular case of Theorem~\ref{gaussible:Usigma_bdd_implies_S_bdd} that $U$ majorizes $S$ in the sense that, for every rearrangement-invariant space $Y(0,1)$,
\begin{equation*}
  \|Sf^*\|_{Y(0,1)} \lesssim \|Uf\|_{Y(0,1)} \quad \text{for every $f$.}
\end{equation*}
Moreover, unlikely $U$, $S$ is linear. The reduction principle leads to a surprising discovery: while the operator $S$, hence the Gaussian--Sobolev embedding, always provides a gain in integrability for example when $X(\mathbb R^n,\gamma_n)=L^p(\mathbb R^n,\gamma_n)$ for $p<\infty$, there is actually a \textit{loss} of integrability for example when $X(\mathbb R^n,\gamma_n)=L^{\infty}(\mathbb R^n,\gamma_n)$ or $X(\mathbb R^n,\gamma_n)=\exp L^{\beta}(\mathbb R^n,\gamma_n)$ for $\beta>0$, in which $\exp L^{\beta}(\mathbb R^n,\gamma_n)$ is the classical exponential-type Orlicz space. Roughly speaking, the gain in integrability vanishes, and eventually is even turned to a loss when we near the endpoint $L^{\infty}$. This is very graphically reflected by the second part of~\eqref{E:endpoints-of-S-restricted}.

In this paper we focus on operators with endpoint behavior similar to that of~\eqref{E:endpoints-of-S} or \eqref{E:endpoints-of-S-lz}, but for a considerably wider class of function spaces of which the operator governing the Gaussian--Sobolev embeddings is a particular instance. Let us finally add for the sake of completeness that operators of another type of nonstandard behavior were studied by different methods in \cite{GP:09,Ma:12}. However, both the motivation of the research and techniques used in those papers were completely different.

\section{Preliminaries}
\begin{conventions}\
        \begin{itemize}
                \item Throughout this paper, $(R, \mu)$ and $(S,\nu)$ are two (possibly different) probabilistic nonatomic measure spaces. If $(R,\mu)=((0,1),\lambda)$, where $\lambda$ is the $1$-dimensional Lebesgue measure on $(0,1)$, we write in short $(0,1)$ instead of $((0,1),\lambda)$.
								\item We write $P\lesssim Q$, where $P,Q$ are nonnegative quantities, when there is a positive constant $c$ independent of all appropriate quantities appearing in the expressions $P$ and $Q$ such that $P\leq c\cdot Q$. If not stated explicitly,  what ``the appropriate quantities appearing in the expressions $P$ and $Q$'' are should be obvious from the context. At the few places where it is not obvious, we will explicitly specify what the appropriate quantities are. We also write $P\gtrsim Q$ with the obvious meaning. Furthermore, we write $P\approx Q$ when $P\lesssim Q$ and $P\gtrsim Q$ simultaneously.
								\item We adhere to the convention that $\frac1{\infty}=0\cdot\infty=0$.
        \end{itemize}
\end{conventions}

We set
\begin{align*}
\M(R, \mu)&= \{f\colon \text{$f$ is a $\mu$-measurable complex-valued function on $R$}\},\\
\intertext{and}
\Mpl(R, \mu)&= \{f \in \M(R, \mu)\colon f \geq 0\ \text{$\mu$-a.e.}\}.
\end{align*}

\subsection*{Rearrangements and rearrangement\hyp{}invariant function spaces}
The \emph{nonincreasing rearrangement} $f^*\colon(0,1) \to [0, \infty]$ of a function $f\in \M(R, \mu)$  is defined as
\begin{equation*}
f^*(t)=\inf\{\lambda\in(0,\infty)\colon\mu\left(\{x\in R\colon|f(x)|>\lambda\}\leq t\}\right),\ t\in(0,1).
\end{equation*}

The \emph{maximal nonincreasing rearrangement} $f^{**} \colon (0,1) \to [0, \infty]$ of a function $f\in \M(R, \mu)$  is
defined as
\begin{equation*}
f^{**}(t)=\frac1t\int_0^ t f^{*}(s)\d{s},\ t\in(0,1).
\end{equation*}
If there is any possibility of misinterpretation, we use the more explicit notations $f^*_\mu$ and $f^{**}_\mu$ instead of $f^*$ and $f^{**}$, respectively, to stress what measure the rearrangements are taken with respect to. The mapping $f \mapsto f^*$ is monotone in the sense that, for every $f,g\in\M(R,\mu)$,
\begin{equation*}
|f| \leq |g| \quad \text{$\mu$-a.e.~on $R$} \quad \Longrightarrow \quad f^* \leq g^* \quad  \text{on $(0, 1)$};
\end{equation*}
consequently, the same implication remains true if ${}^*$ is replaced by ${}^{**}$. We have that $f^*\leq f^{**}$ for every $f\in\M(R,\mu)$.

The \emph{Hardy lemma} (\citep[Chapter~2, Proposition~3.6]{BS}) ensures that, for every $f,g\in\Mpl(0,1)$ and every nonincreasing $h\in\Mpl(0,1)$,
\begin{equation}\label{ch1:ri:hardy-lemma}
\begin{aligned}
&\text{if}\quad\int_0^tf(s)\d{s}\leq \int_0^tg(s)\d{s}\quad\text{for every $t\in(0,1)$,}\\
&\text{then}\quad\int_0^1 f(t)h(t)\d{t}\leq \int_0^1 g(t)h(t)\d{t}.
\end{aligned}
\end{equation}

A functional $\|\cdot\|_{X(0,1)}\colon\Mpl(0,1)\to[0,\infty]$ is called a \emph{rearrangement\hyp{}invariant Banach function norm} (on $(0,1)$) if, for all $f$, $g$ and $\{f_k\}_{k\in \N}$ in $\Mpl(0,1)$, and every $\lambda\in[0,\infty)$:
\begin{itemize}
\item[(P1)] $\|f\|_{X(0,1)}=0$ if and only if $f=0$ a.e.~on $(0,1)$; $\|\lambda f\|_{X(0,1)}= \lambda \|f\|_{X(0,1)}$;  $\|f+g\|_{X(0,1)}\leq \|f\|_{X(0,1)} + \|g\|_{X(0,1)}$;
\item[(P2)] $\|f\|_{X(0,1)}\leq\|g\|_{X(0,1)}$ if $ f\leq g$ a.e.~on $(0,1)$;
\item[(P3)] $\|f_k\|_{X(0,1)} \nearrow \|f\|_{X(0,1)}$ if $f_k \nearrow f$ a.e.~on $(0,1)$;
\item[(P4)] $\|1\|_{X(0,1)}<\infty$;
\item[(P5)] there is a positive constant $C_X$, possibly depending on $\|\cdot\|_{X(0,1)}$ but not on $f$, such that $\int_0^1 f(t)\d{t} \leq C_X \|f\|_{X(0,1)}$;
\item[(P6)] $\|f\|_{X(0,1)} = \|g\|_{X(0,1)}$ whenever $f^*= g^*$.
\end{itemize}

The \textit{Hardy--Littlewood--P\'olya principle} (\citep[Chapter~2, Theorem~4.6]{BS}) asserts that, for every $f,g\in\M(0,1)$ and every rearrangement\hyp{}invariant Banach function norm $\|\cdot\|_{X(0,1)}$,
\begin{equation}\label{ch1:ri:HLP}
\text{if $\int_0^tf^*(s)\d{s}\leq \int_0^tg^*(s)\d{s}$ for every $t\in(0,1)$, then $\|f\|_{X(0,1)}\leq \|g\|_{X(0,1)}$}.
\end{equation}

With every rearrangement\hyp{}invariant Banach function norm $\|\cdot\|_{X(0,1)}$, we associate another functional $\|\cdot\|_{X'(0,1)}$ defined as
\begin{equation*}
\|f\|_{X'(0,1)}= \sup\limits_{\substack{g\in{\Mpl(0,1)}\\\|g\|_{X(0,1)}\leq1}}\int_0^1f(t)g(t)\d{t},\ f\in\Mpl(0,1).
\end{equation*}
The functional $\|\cdot\|_{X'(0,1)}$ is also a rearrangement\hyp{}invariant Banach function norm (\citep[Chapter~2, Proposition~4.2]{BS}), and it is called the \emph{associate Banach function norm} of $\|\cdot\|_{X(0,1)}$. Furthermore, we always have that (\citep[Chapter~1, Theorem~2.7]{BS})
\begin{equation}\label{ch1:ri:normX''}
\|f\|_{X(0,1)}= \sup\limits_{\substack{g\in{\Mpl(0,1)}\\\|g\|_{X'(0,1)}\leq1}}\int_0^1f(t)g(t)\d{t} \quad\text{for every $f\in\Mpl(0,1)$},
\end{equation}
that is,
\begin{equation}\label{ch1:ri:X''=X}
\|\cdot \|_{(X')'(0,1)} = \|\cdot \|_{X(0,1)}.
\end{equation}
The supremum in \eqref{ch1:ri:normX''} does not change when the functions involved are replaced with their nonincreasing rearrangements (\citep[Chapter~2, Proposition~4.2]{BS}), that is,
\begin{equation*}
\|f\|_{X(0,1)}= \sup\limits_{\substack{g\in{\Mpl(0,1)}\\\|g\|_{X'(0,1)}\leq1}}\int_0^1f^*(t)g^*(t)\d{t} \quad\text{for every $f\in\Mpl(0,1)$}.
\end{equation*}

Given a rearrangement\hyp{}invariant Banach function norm $\|\cdot\|_{X(0,1)}$, we define the functional $\|\cdot\|_{X(R,\mu)}$ as
\begin{equation}\label{ch1:ri:rinorm}
\|f\|_{X(R,\mu)}=\|f^*_\mu\|_{X(0,1)} \quad \text{for every $f\in\M(R,\mu)$}.
\end{equation}
Note that $\|f\|_{X(R,\mu)}=\||f|\|_{X(R,\mu)}$. When $(R,\mu)=(0,1)$, \eqref{ch1:ri:rinorm} extends the given rearrangement\hyp{}invariant Banach function norm to all $f\in\M(0,1)$. The functional $\|\cdot\|_{X(R,\mu)}$ restricted to the linear set $X(R,\mu)$ defined as
\begin{equation}\label{ch1:ri:rispace}
X(R,\mu)=\{f\in\M(R,\mu)\colon \|f\|_{X(R,\mu)}<\infty\}
\end{equation}
is a norm (provided that we identify any two functions from $\M(R,\mu)$ coinciding $\mu$-a.e.~on $R$, as usual). In fact, $X(R,\mu)$ endowed with the norm $\|\cdot\|_{X(R,\mu)}$ is a Banach space (\citep[Chapter~1, Theorem~1.6]{BS}). We say that $X(R,\mu)$ is a \emph{rearrangement\hyp{}invariant Banach function space} (an \emph{r.i.~Banach function space}). Note that $f\in\M(R,\mu)$ belongs to $X(R,\mu)$ if and only if $\|f\|_{X(R,\mu)}<\infty$.

The rearrangement\hyp{}invariant Banach function space $X'(R,\mu)$ built upon the associate Banach function norm $\|\cdot\|_{X'(0,1)}$ of $\|\cdot\|_{X(0,1)}$ is called the \emph{associate Banach function space} of $X(R,\mu)$. Thanks to \eqref{ch1:ri:X''=X}, we have that $(X')'(R,\mu)=X(R,\mu)$. Furthermore, one has that
\begin{equation}\label{ch1:ri:holder}
\int_R|f||g|\d{\mu}\leq\|f\|_{X(R,\mu)}\|g\|_{X'(R,\mu)}\quad\text{for every $f,g\in\M(R,\mu)$}.
\end{equation}
We shall refer to \eqref{ch1:ri:holder} as the H\"older inequality.

A functional $\|\cdot\|_{X(0,1)}\colon\Mpl(0,1)\to[0,\infty]$ is called a \emph{rearrangement\hyp{}invariant quasi-Banach function norm} (on $(0,1)$) if it satisfies all the properties of a rearrangement\hyp{}invariant Banach function norm but (P1) and (P5), and instead of (P1) it satisfies, for every $f,g\in\Mpl(0,1)$ and $\lambda\geq0$,
\begin{itemize}
\item[(P1')] $\|f\|_{X(0,1)}=0$ if and only if $f=0$ a.e.~on $(0,1)$; $\|\lambda f\|_{X(0,1)}= \lambda \|f\|_{X(0,1)}$;  there is a constant $C \geq 1$ such that $\|f+g\|_{X(0,1)}\leq C\left(\|f\|_{X(0,1)} + \|g\|_{X(0,1)}\right)$.
\end{itemize}
Given a rearrangement\hyp{}invariant quasi-Banach function norm $\|\cdot\|_{X(0,1)}$, the functional defined by \eqref{ch1:ri:rinorm} is a quasinorm on the linear set defined by \eqref{ch1:ri:rispace}. Moreover, $X(R,\mu)$ endowed with the quasinorm $\|\cdot\|_{X(R,\mu)}$ is a quasi-Banach space (\cite[Corollary~3.7]{PeNe:quasi}), and we called it a \emph{rearrangement\hyp{}invariant quasi-Banach function space} (an \emph{r.i.~quasi-Banach function space}). The rearrangement\hyp{}invariant (quasi-)Banach function space $X(0,1)$ is called the \emph{representation space} of $X(R,\mu)$.

Statements like, ``let $X(R,\mu)$ be a rearrangement\hyp{}invariant (quasi-)Banach function space'', are to be interpreted as ``let $\|\cdot\|_{X(0,1)}$ be a rearrangement\hyp{}invariant (quasi-)Banach function norm and let $X(R,\mu)$ be the corresponding rearrangement\hyp{}invariant (quasi-)Banach function space''.

Let $X(R,\mu)$ and $Y(R,\mu)$ be rearrangement\hyp{}invariant (quasi-)Banach function spaces over the same measure space. We say that $X(R,\mu)$ is \emph{embedded in} $Y(R,\mu)$, and we write $X(R,\mu)\hookrightarrow Y(R,\mu)$, if there is a positive constant $C$ such that $\|f\|_{Y(R,\mu)}\leq C\|f\|_{X(R,\mu)}$ for every $f\in\M(R,\mu)$. If $X(R,\mu)\hookrightarrow Y(R,\mu)$ and $Y(R,\mu)\hookrightarrow X(R,\mu)$ simultaneously, we write $X(R,\mu)=Y(R,\mu)$. We have that (\citep[Chapter~1, Theorem~1.8]{BS} and \cite[Corollary~3.9]{PeNe:quasi})
\begin{equation*}
X(R,\mu)\hookrightarrow Y(R,\mu)\quad\text{if and only if}\quad X(R,\mu)\subseteq Y(R,\mu).
\end{equation*}

We say that a rearrangement\hyp{}invariant quasi-Banach function norm $\|\cdot\|_{X(0,1)}$ is \emph{$p$-convex}, where $p\in(0,\infty)$, if the functional
\begin{equation*}
\|f\|_{X^\frac1{p}(0,1)} = \|f^\frac1{p}\|^p_{X(0,1)},\ f\in\Mpl(0,1),
\end{equation*}
is a rearrangement\hyp{}invariant Banach function norm. The corresponding rearrangement\hyp{}invariant Banach function space $X^\frac1{p}(R,\mu)$ is said to be $p$-convex.

The \emph{$K$-functional} for a couple of (quasi-)Banach function spaces $(X_0(R,\mu), X_1(R,\mu))$ is defined, for every $f \in \M(R,\mu)$ and $t\in(0,\infty)$, as
\begin{equation*}
K(f,t;X_0,X_1)  = \inf\{\lVert g\rVert_{X_0} + t\lVert h\rVert_{X_1}\colon f = g + h\},
\end{equation*}
where the infimum is taken over all representations $f = g + h$ with $g \in X_0$ and $h \in X_1$. If $f\notin (X_0 + X_1)(R,\mu)$, then the infimum is to be interpreted as $\infty$.

\subsection*{Orlicz--Karamata spaces}
A measurable function $\b\colon (0,1)\to(0,\infty)$ is said to be \emph{slowly varying} if for every $\varepsilon > 0$ there are a nondecreasing function $\b_\varepsilon$ and a nonincreasing function $\b_{-\varepsilon}$ such that $t^\varepsilon \b(t) \approx \b_\varepsilon(t)$ and $t^{-\varepsilon} \b(t) \approx \b_{-\varepsilon}(t)$ on $(0,1)$. A slowly varying function $\b$ satisfies
\begin{equation*}
0<\inf_{t\in[a,1)}\b(t)\leq\sup_{t\in[a,1)}\b(t)<\infty \quad \text{for every $a\in(0,1)$}.
\end{equation*}
A positive linear combination of slowly varying functions is a slowly varying function. If $\b_1,\b_2$ are slowly varying functions, so is  $\b_1\b_2$. Any real power of a slowly varying function is a slowly varying function. For every $\alpha\in\R\setminus\{0\}$ and a slowly varying function $\b$, we have that $\lim_{t\to0^+}t^\alpha \b(t)=\lim_{t\to0^+}t^\alpha$. Furthermore, if $\alpha>0$, then
\begin{equation*}
\int_0^t s^{-1+\alpha} \b(s)\d{s} \approx t^\alpha \b(t) \quad \text{for every $t\in (0,1)$}.
\end{equation*}
For more details, we refer the reader to \cite{GoOpTr:05, Pe:21}.

The Orlicz--Karamata functional $\|\cdot\|_{L^{p,\b}(0,1)}$, where $p\in(0,\infty]$ and $\b$ is a slowly varying function, is defined as
\begin{equation*}
\|f\|_{L^{p,\b}(0,1)}=\|\b(t)f^*(t)\|_{L^p(0,1)},\ f\in\Mpl(0,1),
\end{equation*}
where $\|\cdot\|_{L^p(0,1)}$ is the Lebesgue quasi-norm on $(0,1)$, that is,
\begin{equation*}
\|f\|_{L^p(0,1)} = \begin{cases}

\left(\int_0^1|f(t)|^p\d{t}\right)^\frac1{p} \quad &\text{if $p\in(0, \infty)$},\\
\esssup\limits_{t\in(0,1)}|f(t)| \quad &\text{if $p=\infty$}.

\end{cases}
\end{equation*}
The functional $\|\cdot\|_{L^{p,\b}(0,1)}$ is a rearrangement\hyp{}invariant quasi-Banach function norm provided that either $p\in(0,\infty)$ or $p=\infty$ and $\b\in L^\infty(0,1)$ (\cite[Proposition~3.7]{Pe:21}). The corresponding function spaces are called \emph{Orlicz--Karamata spaces}. The Orlicz--Karamata functional $\|\cdot\|_{L^{p,\b}(0,1)}$ is equivalent to a rearrangement\hyp{}invariant Banach function norm if and only if $p=1$ and $\b$ is equivalent to a nonincreasing function, or $p\in(1,\infty)$, or $p=\infty$ and $\b\in L^\infty(0,1)$ (\cite[Theorem~3.26]{Pe:21}). The class of Orlicz--Karamata spaces contains Lebesgue spaces as well as some important Orlicz spaces. If $\b\equiv1$, then $\|\cdot\|_{L^{p,\b}(0,1)}=\|\cdot\|_{L^{p}(0,1)}$ (\cite[Chapter~2, Proposition~1.8]{BS}). Set $\ell(t)=1-\log(t)$, $t\in(0,1)$. If $p\in[1,\infty)$ and $\b=\ell^\alpha$, where $\alpha > 0$ if $p=1$, otherwise $\alpha\in\R$, then $L^{p,\b}(R,\mu)=L^p\left(\log L\right)^{\alpha p}(R,\mu)$, the Orlicz space induced by a Young function $\Phi$ satisfying, for large values of $t$, $\Phi(t)\approx t^p\ell^{\alpha p}(t)$. Furthermore, if $\b=\ell^\alpha$, where $\alpha < 0$, then $L^{\infty,\b}(R,\mu)=\exp L^{-\frac1{\alpha}}(R,\mu)$, the Orlicz space induced by a Young function $\Phi$ satisfying, for large values of $t$, $\Phi(t)\approx \exp(t^{-\frac1{\alpha}})$. For more details, we refer the reader to \cite[Section~8]{OpPi:99}.

\section{Main results}

In this section we shall state and proof our main results. We begin by introducing a key function.

\begin{definition}\label{D:sigma}
\noindent Let $p\in(0, \infty)$ and $\b_1, \b_2$ be slowly varying functions. We define the function
\begin{equation*}
\sigma=\sigma(\b_1, \b_2,p)\colon[0,1]\to[0,1]
\end{equation*}
as the increasing, bijective function satisfying
\begin{equation}\label{sigma_integral_identity}
    t^p = \frac 1C\int_0^{\sigma(t)^p}\left[\b_1(s) \b_2(s)^{-1}\right]^p\d{s}\quad\text{for every $t\in[0,1]$,}
\end{equation}
where
\begin{equation*}
C = \int_0^1 \left[\b_1(s)\b_2(s)^{-1}\right]^p\d{s}  \in (0,\infty).
\end{equation*}
\end{definition}

\nopagebreak\begin{remark}
It immediately follows from \eqref{sigma_integral_identity} that $\sigma,\sigma^{-1}\in\mathcal C([0,1])$, where $\sigma^{-1}$ denotes the inverse function. Furthermore, we have that
\begin{equation*}
\quad \sigma^{-1}(t) \approx t\b_1(t^p)\b_2(t^p)^{-1}
\end{equation*}
on $(0,1)$. If the functions $\b_1, \b_2$ are continuous, then $\sigma,\sigma^{-1}\in\mathcal C^1(0,1)$ and
\begin{equation*}
\big(\sigma^{-1}(t^{1/p})^p\big)' \approx  \b_1(t)^p \b_2(t)^{-p}
\end{equation*}
on $(0,1)$. We shall use these properties of $\sigma$ without making any explicit reference to them.
\end{remark}

We shall now characterize the $K$-inequality corresponding to the couples $(L^{p,\b_1}, L^{\infty, \b_2})$ and $(L^p,L^\infty)$ by an inequality for certain integrals.

Note that, since $(R,\mu)$ and $(S,\nu)$ are finite nonatomic measure spaces, we have the embeddings $L^{\infty, \b_2}(S,\nu)\hookrightarrow L^{p,\b_1}(S,\nu)$ and $L^\infty(R,\mu)\hookrightarrow L^p(R,\mu)$, and so to write $f\in L^p(R,\mu)$ and $g\in L^{p,\b_1}(S,\nu)$ is the same as to write $f\in \left(L^p+L^\infty\right)(R,\mu)$ and $g\in \left(L^{p,\b_1} + L^{\infty,\b_2}\right)(S,\nu)$.

\begin{theorem}\label{T:gaussible:Kinequality_equiv_to_integral_one}
Let $p\in(0, \infty)$ and $(\b_1,\b_2)$ be a pair of continuous slowly varying functions. Assume that $\b_2$ is nondecreasing. Let $f\in L^p(R,\mu)$ and $g\in L^{p,\b_1}(S,\nu)$. The inequality
\begin{equation}\label{gaussible:Kinequality_equiv_to_integral_one:K_ineq}
K(g, t; L^{p,\b_1}, L^{\infty, \b_2}) \lesssim K(f,t;L^p,L^\infty)
\end{equation}
holds for every $t\in(0,1)$ with a multiplicative constant independent of $f$ and $g$ if and only if the inequality
\begin{equation*}
\int_0^t \left[g^*(s) \b_1(s) \right]^p \,ds \lesssim \int_0^t \left[f^*\big(\sigma^{-1}(s^{1/p})^p\big) \b_1(s) \b_2(s)^{-1}\right]^p \, ds
\end{equation*}
holds for every $t\in(0,1)$ with a multiplicative constant independent of $f$ and $g$.
\end{theorem}

\begin{proof}
First, we claim that, for every $g \in L^{p,\b_1}(S,\nu)$ and every $t\in(0,1)$,
\begin{equation}\label{gaussible:K_functional_sigma_lemma:equivalence}
K(g, t; L^{p,\b_1}, L^{\infty, \b_2}) \approx \left(\int_0^{\sigma(t)^p} \left[g^*(s) \b_1(s) \right]^p\,ds\right)^\frac 1p + t\cdot \sup_{\sigma(t)^p \leq s < 1} g^*(s)\b_2(s),
\end{equation}
in which the multiplicative constants are independent of $g$ and $t$. This result can be derived from a general implicit formula appearing in~\cite[Theorem~4.1]{Mi:79}. Since we need the explicit formula here, we shall prove it in detail for the sake of completeness. Let $g\in L^{p,\b_1}(S,\nu)$ and $t\in(0,1)$. For the sake of brevity, we set
\begin{equation*}
I(g)(t) = \left(\int_0^{\sigma(t)^p} \left[g^*(s) \b_1(s) \right]^p\,ds\right)^\frac 1p + t\cdot \sup_{\sigma(t)^p \leq s < 1} g^*(s)\b_2(s).
\end{equation*}

Let $g = g_1 + g_2$, where $g_1 \in L^{p,\b_1}(S,\nu)$ and $g_2 \in L^{\infty,\b_2}(S,\nu)$, be a decomposition of $g$. We have that
\begin{align}
I(g)(t) &= \left(\int_0^{\sigma(t)^p} \left[(g_1 + g_2)^*(s) \b_1(s)\right]^p\,ds\right)^\frac 1p + t\cdot \sup_{\sigma(t)^p \leq s < 1} (g_1 + g_2)^*(s)\b_2(s) \notag \\
\begin{split}\label{gaussible:K_functional_sigma_lemma:eq1}
		& \leq \left(\int_0^{\sigma(t)^p} \left[(g_1^*(s/2) + g_2^*(s/2)) \b_1(s)\right]^p\,ds\right)^\frac 1p + t\cdot \sup_{\sigma(t)^p \leq s < 1} \left[g_1^*(s/2) + g_2^*(s/2)\right]\b_2(s) \\
		& \lesssim I(g_1^*(\cdot/2))(t) + I(g_2^*(\cdot/2))(t).
\end{split}
\end{align}
As for $I(g_1^*(\cdot/2))(t)$, since
\begin{equation*}
g_1^*(s)\approx g_1^*(s)\left(\frac1{s\b_1(s)^p}\int_0^s  \b_1(\tau)^p\,d\tau \right)^\frac 1p\leq \frac 1{s^\frac 1p\b_1(s)}\left(\int_0^s \left[g_1^*(\tau) \b_1(\tau)\right]^p\,d\tau\right)^\frac 1p
\end{equation*}
for every $s\in(0,1)$, we obtain that
\begin{equation}\label{gaussible:K_functional_sigma_lemma:eq2}
\begin{aligned}
    I(g_1^*(\cdot/2))(t) &\lesssim \lVert g_1 \rVert_{L^{p,\b_1}(S,\nu)} + t\cdot \sup_{\sigma(t)^p \leq s < 1} g_1^*(s/2)\b_2(s) \\
		&\lesssim \lVert g_1 \rVert_{L^{p,\b_1}(S,\nu)} + t\cdot \sup_{\sigma(t)^p \leq s < 1} \frac{\displaystyle \left(\int_0^\frac{s}{2}\left[g_1^*(\tau)\b_1(\tau)\right]^p\,d\tau\right)^\frac 1p}{s^\frac 1p\b_1(s)\b_2(s)^{-1}}\\
		&\leq \lVert g_1 \rVert_{L^{p,\b_1}(S,\nu)}\left(1+ t\cdot \sup_{\sigma(t)^p \leq s < 1} \frac1{s^\frac 1p\b_1(s)\b_2(s)^{-1}}\right)\\ &\approx \lVert g_1 \rVert_{L^{p,\b_1}(S,\nu)}\left(1+ t\cdot \frac1{\sigma(t)\b_1(\sigma(t)^p)\b_2(\sigma(t)^p)^{-1}}\right)
		\\ &\approx \lVert g_1 \rVert_{L^{p,\b_1}(S,\nu)}.
\end{aligned}
\end{equation}
As for the second term on the right-hand side of \eqref{gaussible:K_functional_sigma_lemma:eq1}, note that
\begin{equation}\label{gaussible:K_functional_sigma_lemma:eq3}
\begin{aligned}
    I(g_2^*(\cdot/2))(t) &\lesssim \lVert g_2 \rVert_{L^{\infty, \b_2}(S,\nu)} \left(\int_0^{\sigma(t)^p}\left[\b_1(s) \b_2(s)^{-1}\right]^p\,ds\right)^\frac 1p + t\cdot \sup_{\sigma(t)^p \leq s < 1} g_2^*(s/2) \b_2(s/2) \\
		&\approx t\lVert g_2 \rVert_{L^{\infty, \b_2}(S,\nu)} + t\cdot \sup_{\frac{\sigma(t)^p}{2} \leq s < \frac{1}{2}} g_2^*(s) \b_2(s) \\
		&\lesssim t\cdot \lVert g_2 \rVert_{L^{\infty, \b_2}(S,\nu)},
\end{aligned}
\end{equation}
in which we used \eqref{sigma_integral_identity}. Hence, by combining \eqref{gaussible:K_functional_sigma_lemma:eq3} and \eqref{gaussible:K_functional_sigma_lemma:eq2} together with \eqref{gaussible:K_functional_sigma_lemma:eq1}, and taking the infimum over all such representations $g = g_1 + g_2$, we obtain that
\begin{equation}\label{E:upper-bound}
I(g)(t) \lesssim K(g, t; L^{p,\b_1}, L^{\infty,\b_2}).
\end{equation}

As for the opposite inequality, we may assume that $I(g)(t)<\infty$, for otherwise there is nothing to prove. Define the functions $g_1,g_2\in\M(S,\nu)$ as
\begin{equation*}
g_1(x) = \max\{|g(x)| - g^*(\sigma(t)^p), 0\}\cdot \sgn g(x),\ x\in S,
\end{equation*}
and
\begin{equation*}
g_2(x) = g(x) - g_1(x) = \min\{|g(x)|, g^*(\sigma(t)^p)\}\cdot \sgn g(x),\ x\in S.
\end{equation*}
Clearly, $g=g_1+g_2$ and we have that
\begin{equation*}
g_1^*(s) = (g^*(s) - g^*(\sigma(t)^p))\chi_{(0,\sigma(t)^p)}(s) \quad \text{ and } \quad g_2^*(s) = \min\{g^*(s), g^*(\sigma(t)^p)\}
\end{equation*}
for every $s\in(0,1)$. Note that
\begin{equation*}
    \lVert g_1 \rVert_{L^{p,\b_1}(S,\nu)} = \left(\int_0^{\sigma(t)^p} \left[(g^*(s) - g^*(\sigma(t)^p))\b_1(s)\right]^p\, ds\right)^\frac 1p \leq \left(\int_0^{\sigma(t)^p} \left[g^*(s)\b_1(s)\right]^p\, ds\right)^\frac 1p
\end{equation*}
and
\begin{align*}
    \lVert g_2 \rVert_{L^{\infty, \b_2}(S,\nu)} &= g^*(\sigma(t)^p) \sup_{0 < s \leq \sigma(t)^p} \b_2(s) + \sup_{\sigma(t)^p \leq s < 1} g^*(s)\b_2(s) \approx \sup_{\sigma(t)^p \leq s < 1} g^*(s)\b_2(s);
\end{align*}
consequently,
\begin{equation*}
\lVert g_1 \rVert_{L^{p,\b_1}(S,\nu)} + t\lVert g_2 \rVert_{L^{\infty, \b_2}(S,\nu)} \lesssim I(g)(t)
\end{equation*}
and $g_1\in L^{p,\b_1}(S,\nu)$ and $g_2\in L^{\infty, \b_2}(S,\nu)$. Hence,
\begin{equation*}
K(g, t; L^{p,\b_1},L^{\infty, \b_2}) \lesssim I(g)(t),
\end{equation*}
which together with~\eqref{E:upper-bound} establishes our claim~\eqref{gaussible:K_functional_sigma_lemma:equivalence}.

Second, since we have that (see~\cite{Kree:67}, also~\cite[Theorem~4.1]{Holmstedt:70})
\begin{equation*}
K(f,t;L^p,L^\infty)\approx\left(\int_0^{t^p}f^*(s)^p \,ds\right)^\frac 1p,
\end{equation*}
in which the multiplicative constants depend only on $p$, in view of~\eqref{gaussible:K_functional_sigma_lemma:equivalence}, we need to prove that
\begin{equation}\label{gaussible:Kinequality_equiv_to_integral_one:eq1}
\left(\int_0^{\sigma(t)^p} \left[g^*(s) \b_1(s) \right]^p\,ds\right)^\frac 1p + t\cdot \sup_{\sigma(t)^p \leq s < 1} g^*(s)\b_2(s) \lesssim \left(\int_0^{t^p}f^*(s)^p\,ds\right)^\frac 1p
\end{equation}
for every $t\in(0,1)$ if and only if
\begin{equation}\label{gaussible:Kinequality_equiv_to_integral_one:eq2}
\int_0^t \left[g^*(s) \b_1(s) \right]^p \,ds \lesssim \int_0^t \left[f^*\big(\sigma^{-1}(s^{1/p})^p\big) \b_1(s) \b_2(s)^{-1}\right]^p \, ds \quad \text{for every $t\in(0,1)$.}
\end{equation}
 We shall observe that
\begin{equation}\label{gaussible:Kinequality_equiv_to_integral_one:eq3}
\int_0^{\sigma(t)^p} \left[g^*(s) \b_1(s) \right]^p\,ds \lesssim \int_0^{t^p}f^*(s)^p\,ds \quad \text{for every $t\in(0,1)$}
\end{equation}
if and only if
\begin{equation}\label{gaussible:Kinequality_equiv_to_integral_one:eq4}
t^p\cdot\sup_{\sigma(t)^p \leq s < 1} \frac1{s\b_1(s)^p\b_2(s)^{-p}}\int_0^s \left[g^*(\tau)\b_1(\tau)\right]^p\,d\tau \lesssim \int_0^{t^p} f^*(s)^p\,ds \quad \text{for every $t\in(0,1)$}.
\end{equation}
Note that \eqref{gaussible:Kinequality_equiv_to_integral_one:eq4} plainly implies \eqref{gaussible:Kinequality_equiv_to_integral_one:eq3} inasmuch as $$\sigma(t)\b_1(\sigma(t)^p)\b_2(\sigma(t)^p)^{-1}\approx t \quad \text{for every $t\in(0,1)$}.$$ As for the opposite implication, let $t\in(0,1)$. Since the function $(0,1)\ni t\mapsto (|f|^p)^{**}(t^p)$ is nonincreasing, \eqref{gaussible:Kinequality_equiv_to_integral_one:eq3} actually implies that
\begin{equation}\label{gaussible:Kinequality_equiv_to_integral_one:eq5}
\sup_{t \leq s < 1} \frac1{s^p}\int_0^{\sigma(s)^p} \left[g^*(\tau)\b_1(\tau)\right]^p\,d\tau \lesssim \frac1{t^p}\int_0^{t^p}f^*(s)^p\,ds.
\end{equation}
Since $\sigma^{-1}$ is an increasing bijection of $[0,1]$ onto itself, by the change of variables $s = \sigma^{-1}(\tilde s^{1/p})$, \eqref{gaussible:Kinequality_equiv_to_integral_one:eq5} is equivalent to
\begin{equation*}
\sup_{\sigma(t)^p \leq s < 1} \frac1{s\b_1(s)^p\b_2(s)^{-p}}\int_0^s \left[g^*(\tau)\b_1(\tau)\right]^p\,d\tau \lesssim \frac1{t^p}\int_0^{t^p}f^*(s)^p\,ds,
\end{equation*}
whence \eqref{gaussible:Kinequality_equiv_to_integral_one:eq4} follows. Furthermore, by the change of variables $s = \sigma^{-1}(\tilde s^{1/p})^p$, we have that
\begin{equation*}
\int_0^{t^p}f^*(s)\,ds \approx \int_0^{\sigma(t)^p}f^*(\sigma^{-1}(s^{1/p})^p)\b_1(s)^p\b_2(s)^{-p}\,ds \quad \text{for every $t\in(0,1)$}.
\end{equation*}
Hence, since $\sigma$ is a bijection of $[0,1]$ onto itself, \eqref{gaussible:Kinequality_equiv_to_integral_one:eq2} is equivalent to \eqref{gaussible:Kinequality_equiv_to_integral_one:eq3}.

Finally, the proof will be completed once we show that \eqref{gaussible:Kinequality_equiv_to_integral_one:eq4} is equivalent to \eqref{gaussible:Kinequality_equiv_to_integral_one:eq1}. Since \eqref{gaussible:Kinequality_equiv_to_integral_one:eq1} plainly implies \eqref{gaussible:Kinequality_equiv_to_integral_one:eq3}, which is equivalent to \eqref{gaussible:Kinequality_equiv_to_integral_one:eq4}, we only need to observe that \eqref{gaussible:Kinequality_equiv_to_integral_one:eq4} implies \eqref{gaussible:Kinequality_equiv_to_integral_one:eq1} (the former actually implies the latter pointwise). To this end, note that
\begin{align*}
\sup_{\sigma(t)^p \leq s < 1} g^*(s)\b_2(s) &\approx \sup_{\sigma(t)^p \leq s < 1} \frac{g^*(s)}{s^\frac 1p\b_1(s)\b_2(s)^{-1}} \left(\int_0^s \b_1(\tau)^p\,d\tau\right)^\frac 1p \\
& \leq \sup_{\sigma(t)^p \leq s < 1} \frac 1{s^\frac 1p\b_1(s)\b_2(s)^{-1}}\left( \int_0^s \left[g^*(\tau)\b_1(\tau)\right]^p\,d\tau\right)^\frac 1p
\end{align*}
for every $t\in(0,1)$. Hence, if \eqref{gaussible:Kinequality_equiv_to_integral_one:eq4} is true (consequently, so is \eqref{gaussible:Kinequality_equiv_to_integral_one:eq3}), then
\begin{align*}
&K(g, t; L^{p,\b_1}, L^{\infty, \b_2}) \\&\lesssim \left(\int_0^{\sigma(t)^p} \left[g^*(s) \b_1(s) \right]^p\,ds\right)^\frac 1p + t\cdot \sup_{\sigma(t)^p \leq s < 1} \frac 1{s^\frac 1p\b_1(s)\b_2(s)^{-1}}\left( \int_0^s \left[g^*(\tau)\b_1(\tau)\right]^p\,d\tau\right)^\frac 1p \\ &\lesssim \left(\int_0^{t^p}f^*(s)^p\,ds\right)^\frac 1p,
\end{align*}
for every $t\in(0,1)$.

\end{proof}

\begin{remark}
If \eqref{gaussible:Kinequality_equiv_to_integral_one:K_ineq} is valid for every $t\in(0,1)$, it is actually valid for every $t\in(0, \infty)$ (with a possibly different multiplicative constant). Indeed, owing to the embeddings mentioned above Theorem~\ref{T:gaussible:Kinequality_equiv_to_integral_one}, we have that $K(f,t;L^p,L^\infty)\approx K(f,1;L^p,L^\infty)$ and $K(g,t;L^{p,\b_1},L^{\infty, \b_2})\approx K(g,1;L^{p,\b_1}, L^{\infty, \b_2})$ for every $t\in[1,\infty)$, in which the multiplicative constants are independent of $f,g$ and $t$; therefore
\begin{align*}
K(g,t;L^{p,\b_1},L^{\infty, \b_2}) &\approx K(g,1;L^{p,\b_1},L^{\infty, \b_2}) \\
&\lesssim K(f,1;L^p,L^\infty) \approx K(f,t;L^p,L^\infty)
\end{align*}
for every $t\in[1,\infty)$.
\end{remark}

Now we shall introduce a key notion of a gaussible operator.

\begin{definition}
Let $p\in(0, \infty)$ and $\b_1,\b_2$ be slowly varying functions. We say that an operator $T$ defined on $L^p(R,\mu)$ having values in $\M(S,\nu)$ is \emph{$(p,\b_1,\b_2)$-gaussible} if
\begin{equation*}
\int_0^t \left[(Tf)^*(s)\b_1(s)\right]^p \, ds \lesssim \int_0^t \left[f^*\big(\sigma^{-1}(s^{1/p})^p\big) \b_1(s) \b_2(s)^{-1}\right]^p \, ds
\end{equation*}
\noindent for every $f\in L^p(R,\mu)$ and $t\in(0,1)$.
\end{definition}

\begin{remarks}\label{gaussible:gaussible_maps_to_Lpb1}\hphantom{}
\begin{enumerate}[(i)]
\item It follows immediately from the definition that a $(p,\b_1,\b_2)$-gaussible operator is bounded from $L^p(R,\mu)$ to $L^{p,\b_1}(S,\nu)$. Indeed, any $(p,\b_1,\b_2)$-gaussible operator $T$ satisfies
\begin{equation*}
\|Tf\|_{L^{p,\b_1}(S,\nu)} \lesssim \left(\int_0^1\left[f^*\big(\sigma^{-1}(s^{1/p})^p\big)\b_1(s)\b_2(s)^{-1}\right]^p\,ds\right)^{1/p} \approx \|f\|_{L^p(R,\mu)}.
\end{equation*}
\item In view of Theorem~\ref{T:gaussible:Kinequality_equiv_to_integral_one}, an operator $T$ defined on $L^p(R,\mu)$ having values in $\M(S,\nu)$ is $(p,\b_1,\b_2)$-gaussible if and only if it satisfies
\begin{equation}\label{gaussible:gaussible_maps_to_Lpb1:K_ineq}
K(Tf, t; L^{p,\b_1}, L^{\infty, \b_2}) \lesssim K(f,t;L^p,L^\infty) \quad \text{for every $f\in L^p(R,\mu)$ and $t\in (0,1)$}.
\end{equation}
\item The class of operators satisfying the $K$-inequality \eqref{gaussible:gaussible_maps_to_Lpb1:K_ineq} actually coincides with a certain class of operators introduced in \cite[Section~4.1]{BKbook}. An operator $T$ defined on $X_0 + X_1$ having values in $Y_0 + Y_1$, where $(X_0,X_1)$ and $(Y_0,Y_1)$ are two pairs of compatible couples of quasi-Banach spaces, belongs to the class $B(X_0, X_1; Y_0, Y_1)$ if there is a constant $C>0$ such that, for every $f_i\in X_i$, $i=0,1$, and every $\varepsilon > 0$, there are $g_i\in Y_i$, $i=0,1$, such that
\begin{equation*}
T(f_0 + f_1) = g_0 + g_1 \quad \text{and} \quad \|g_i\|_{Y_i}\leq C \|f_i\|_{X_i} + \varepsilon,\ i=0,1.
\end{equation*}
By \cite[Proposition~4.1.3]{BKbook} with some appropriate modifications, an operator $T$ defined on $L^p(R,\mu)$ having values in $\M(S,\nu)$ satisfies the $K$-inequality \eqref{gaussible:gaussible_maps_to_Lpb1:K_ineq} if and only if it belongs to the class $B(L^p(R,\mu), L^\infty(R,\mu); L^{p,\b_1}(S,\nu), L^{\infty, \b_2}(S,\nu))$.

Indeed, assume that \eqref{gaussible:gaussible_maps_to_Lpb1:K_ineq} holds. Let $f_0\in L^p(R,\mu)$, $f_1\in L^\infty(R,\mu)$ and $\varepsilon>0$ be given. Assume that neither $f_0$ nor $f_1$ is equivalent to the zero function (otherwise the proof is trivial). Thanks to \eqref{gaussible:gaussible_maps_to_Lpb1:K_ineq} with $t=t_0=\dfrac{\lVert f_0 \rVert_{L^p(R,\mu)}}{\lVert f_1 \rVert_{L^\infty(R,\mu)}}$, there are $g_0\in L^{p,\b_1}(S,\nu)$ and $g_1\in L^{\infty, \b_2}(S,\nu)$ such that $T(f_0+f_1)=g_0+g_1$ and
\begin{equation*}
\lVert g_0\rVert_{L^{p,\b_1}(S,\nu)} + \dfrac{\lVert f_0 \rVert_{L^p(R,\mu)}}{\lVert f_1 \rVert_{L^\infty(R,\mu)}} \lVert g_1\rVert_{L^{\infty,\b_2}(S,\nu)} \le
2C\lVert f_0 \rVert_{L^p(R,\mu)} + \min\left\{ \dfrac{\lVert f_0 \rVert_{L^p(R,\mu)}}{\lVert f_1 \rVert_{L^\infty(R,\mu)}}, 1 \right\}\varepsilon,
\end{equation*}
whence
\begin{align*}
\lVert g_0\rVert_{L^{p,\b_1}(S,\nu)} &\leq 2C\lVert f_0 \rVert_{L^p(R,\mu)} + \varepsilon\\
\intertext{and}
\lVert g_1\rVert_{L^{\infty,\b_2}(S,\nu)} &\leq 2C\lVert f_1 \rVert_{L^\infty(R,\mu)} + \varepsilon.
\end{align*}
Hence $T\in B(L^p(R,\mu), L^\infty(R,\mu); L^{p,\b_1}(S,\nu), L^{\infty, \b_2}(S,\nu))$. Conversely, assume that $T\in B(L^p(R,\mu), L^\infty(R,\mu); L^{p,\b_1}(S,\nu), L^{\infty, \b_2}(S,\nu))$. Let $f \in L^p(R,\mu)$ and $t\in(0,\infty)$ be given. Let $f=f_0+f_1$ be a decomposition of $f$, where $f_0 \in L^p(R,\mu)$, $f_1 \in L^\infty(R, \mu)$. Fix arbitrary $\varepsilon > 0$. There are $g_0 \in L^{p,\b_1}(S,\nu)$ and $g_1 \in L^{\infty,\b_2}(S,\nu)$ such that $T(f_0+f_1)=g_0+g_1$ and
\begin{equation*}
\lVert g_0\rVert_{L^{p,\b_1}(S,\nu)} \leq C \lVert f_0 \rVert_{L^p(R,\mu)} + \varepsilon \quad \text{ and } \quad \lVert g_1\rVert_{L^{\infty,\b_2}(S,\nu)} \leq C \lVert f_1 \rVert_{L^\infty(R,\mu)} + \varepsilon,
\end{equation*}
where $C>0$ is a constant independent of $f_0$, $f_1$, $g_0$, $g_1$, $t$ and $\varepsilon$. Consequently,
\begin{align*}
K(Tf, t; L^{p,\b_1}, L^{\infty,\b_2}) &\leq
C(\lVert f_0 \rVert_{L^p(R,\mu)} + t\lVert f_1 \rVert_{L^\infty(R,\mu)}) + (1+t)\varepsilon.
\end{align*}
Since $\varepsilon>0$ was arbitrary, it follows that
\begin{equation*}
K(Tf, t; L^{p,\b_1}, L^{\infty,\b_2}) \leq C(\lVert f_0 \rVert_{L^p(R,\mu)} + t\lVert f_1 \rVert_{L^\infty(R,\mu)}).
\end{equation*}
By taking the infimum over all decompositions $f=f_0+f_1$, $f_0 \in L^p(R,\mu)$, $f_1 \in L^\infty(R, \mu)$, we obtain \eqref{gaussible:gaussible_maps_to_Lpb1:K_ineq}.
\end{enumerate}
\end{remarks}

We now specify the class of pairs of slowly varying functions for which we shall later obtain our main result.

\begin{definition}\label{D:bp}
\noindent Let $p\in(0, \infty)$. We say that a pair $(\b_1,\b_2)$ of slowly varying functions belongs to the class $\mathcal B_p$ if

\begin{enumerate}[(a)]
    \item $\b_1, \b_2\in\mathcal C(0,1)$,
		
	\item $\b_1$ is nonincreasing and $\b_2$ is nondecreasing,

    \item $\b_1(t) \approx \b_1\left(t\b_1(t)^p\b_2(t)^{-p}\right)$ near $0^+$,

    \item $\displaystyle \sup_{0 < t < 1} \b_2(t)^p \int_t^1 \frac{ds}{s \b_1(s)^p} < \infty$.
\end{enumerate}

\end{definition}

\begin{remarks}\label{gaussible:rem:properties_of_sigma}\phantom{}
\begin{enumerate}[(i)]

    \item Note that (c) in Definition~\ref{D:bp} actually implies that
\begin{equation*}
\b_1(t) \approx \b_1\big(\sigma^{-1}(t^{1/p})^p\big) \approx \b_1\big(\sigma(t^{1/p})^p\big) \quad \text{for every $t\in(0,1)$}.
\end{equation*}

    \item Since the function $t\mapsto \b_1(t)^p \b_2(t)^{-p}$, $t\in(0,1)$, is nonincreasing, it follows that
\begin{equation}\label{tlessthansigma-1}
t\leq \sigma^{-1}(t^{1/p})^p \quad \text{for every $t\in(0,1)$}.
\end{equation}
Indeed, owing to \eqref{sigma_integral_identity}, we have that
\begin{align*}
\frac{\sigma^{-1}(t^{1/p})^p}{t} &=\left(\int_0^1 \left[\b_1(s)\b_2(s)^{-1}\right]^p\d{s}\right)^{-1}\frac1{t}\int_0^{t}\left[\b_1(s) \b_2(s)^{-1}\right]^p\d{s} \\
&\geq\left(\int_0^1 \left[\b_1(s)\b_2(s)^{-1}\right]^p\d{s}\right)^{-1}\int_0^1\left[\b_1(s) \b_2(s)^{-1}\right]^p\d{s}\\
&=1
\end{align*}
for every $t\in(0,1)$. Moreover, if the function $t\mapsto \b_1(t)^p \b_2(t)^{-p}$, $t\in(0,1)$, is decreasing, then the inequality in \eqref{tlessthansigma-1} is strict.
\end{enumerate}
\end{remarks}

Now we shall introduce three operators, which will play an essential role in what follows.

\begin{definition}
Let $p\in(0, \infty)$ and $\b_1, \b_2$ be slowly varying functions. We define the operators $U_{\b_1, \b_2,p}$, $T_{\b_1, \b_2,p}$ and $S_{\b_1, p}$ as, for every $f\in\M(0,1)$,
\begin{align*}
U_{\b_1, \b_2,p} f(t)&= f^*(\sigma^{-1}(t^{1/p})^p) \b_2(t)^{-1},\ t\in(0,1),\\
T_{\b_1, \b_2, p} f(t)&=\sup_{t\leq s < 1}\frac{f^*(\sigma(s^{1/p})^p)}{\b_1(\sigma(s^{1/p})^p)^p},\ t\in(0,1),\\
\intertext{and}
S_{\b_1, p}f(t)&=\left(\int_t^1\frac{|f(s)|^p}{s\b_1(s)^p}\, ds\right)^\frac 1p,\ t\in(0,1),
\end{align*}
where $\sigma$ is the function from Definition~\ref{D:sigma}.
\end{definition}

\begin{remarks}\phantom{}
\begin{enumerate}[(i)]
\item Note that the functions $U_{\b_1, \b_2,p} f$, $T_{\b_1, \b_2, p} f$ and $S_{\b_1, p}f$ are nonincreasing for every $f\in\M(0,1)$.

\item The operator $U_{\b_1, \b_2,p}$ is plainly $(p,\b_1,\b_2)$-gaussible (with $(R,\mu)=(S,\nu)=(0,1)$). Hence

\begin{equation*}
K(U_{\b_1, \b_2,p}f, t; L^{p,\b_1}, L^{\infty, \b_2}) \lesssim K(f,t;L^p,L^\infty) \quad \text{for every $f\in L^p(R,\mu)$ and $t\in (0,1)$}
\end{equation*}
owing to Remark~\ref{gaussible:gaussible_maps_to_Lpb1}(ii). Moreover, although $U_{\b_1, \b_2,p}$ is neither linear nor quasilinear, it is bounded (in the classic sense) from $L^p(0,1)$ to $L^{p,\b_1}(0,1)$ (see Remark~\ref{gaussible:gaussible_maps_to_Lpb1}(i)) and from $L^\infty(0,1)$ to $L ^{\infty,\b_2}(0,1)$; indeed,
\begin{equation*}
\|U_{\b_1, \b_2,p} f\|_{L^{\infty,\b_2}(0,1)} = \sup_{0 < t < 1}\left(\frac{f^*(\sigma^{-1}(t^{1/p})^p)}{\b_2(t)}\right)\b_2(t) = \|f\|_{L^\infty(0,1)}
\end{equation*}
for every $f\in\M(0,1)$.
\end{enumerate}
\end{remarks}

Now we are in a position to state and prove our main results.

\begin{theorem}\label{gaussible:calderon_thm:restated}

\noindent Let $p\in(0, \infty)$ and $(\b_1, \b_2) \in \mathcal B_p$. Let $X(R,\mu)$ and $Y(S,\nu)$ be r.i.~quasi-Banach function spaces that are $p$-convex. The following four statements are equivalent.

\begin{enumerate}[(i)]
    \item Every $(p,\b_1,\b_2)$-gaussible operator $T$ is bounded from $X(R,\mu)$ to $Y(S,\nu)$.

    \item Every operator $T$ defined on $L^p(R,\mu)$ having values in $\M(S,\nu)$ that satisfies
		\begin{equation*}
			K(Tf, t; L^{p,\b_1}, L^{\infty, \b_2}) \lesssim K(f,t;L^p,L^\infty) \quad \text{for every $f\in L^p(R,\mu)$ and $t\in (0,1)$}
		\end{equation*}
		is bounded from $X(R,\mu)$ to $Y(S,\nu)$.

    \item The operators $U_{\b_1, \b_2,p}$ and $S_{\b_1, p}$ are bounded from $X(0,1)$ to $Y(0,1)$.
				
	\item The operator $T_{\b_1, \b_2, p}$ is bounded from $\big(Y^\frac 1p\big)'(0,1)$ to $\big(X^\frac 1p\big)'(0,1)$.

\end{enumerate}

\end{theorem}

\begin{proof}
\emph{(i) and (ii) are equivalent.} This is an immediate consequence of the very definition of $(p,\b_1,\b_2)$-gaussible operators and Theorem~\ref{T:gaussible:Kinequality_equiv_to_integral_one}, as was already observed in Remark~\ref{gaussible:gaussible_maps_to_Lpb1}(ii).

\emph{(i) implies (iii).} Note that, for every $f\in L^p(R,\mu)$, the function $U_{\b_1, \b_2,p}(f^*)$ is a nonnegative, nonincreasing, finite function on $(0,1)$. By \cite[Chapter~2, Corollary 7.8]{BS}, there is a function $g_f\in\M(S,\nu)$ such that $g_f^*=U_{\b_1, \b_2,p}(f^*)$. The auxiliary operator $T$ defined as $T f= g_f$, $f\in L^p(R,\mu)$, is plainly $(p,\b_1,\b_2)$-gaussible (note that this does not depend on particular choices of $g_f$). Hence, owing to $(i)$, $T$ is bounded from $X(R,\mu)$ to $Y(S,\nu)$. By \cite[Chapter~2, Corollary 7.8]{BS} again, for every $h\in X(0,1)$, there is a function $f_h\in X(R,\mu)$ such that $f_h^*=h^*$. Therefore,
\begin{equation*}
\|U_{\b_1, \b_2,p} h\|_{Y(0,1)} = \|(T f_h)^*\|_{Y(0,1)} = \|T f_h\|_{Y(S,\nu)} \lesssim \|f_h\|_{X(R,\mu)} = \|h^*\|_{X(0,1)} = \|h\|_{X(0,1)}
\end{equation*}
for every $h\in X(0,1)$. Hence $U_{\b_1, \b_2,p}$ is bounded from $X(0,1)$ to $Y(0,1)$.

Next, it is easy to see that $S_{\b_1, p}$ is bounded from $L^p(0,1)$ to $L^{p,\b_1}(0,1)$ and from $L^\infty(0,1)$ to $L^{\infty,\b_2}(0,1)$. Moreover, the operator is quasilinear. It follows that
\begin{equation*}
K(S_{\b_1, p}f, t; L^{p,\b_1}, L^{\infty, \b_2}) \lesssim K(f,t;L^p,L^\infty) \quad \text{for every $f\in L^p(0,1)$ and $t\in (0,1)$};
\end{equation*}
hence $S_{\b_1, p}$ is a $(p,\b_1,\b_2)$-gaussible operator with respect to $(R,\mu)=(S,\nu)=(0,1)$ (see Remark~\ref{gaussible:gaussible_maps_to_Lpb1}(ii)). Arguing along the same lines as for $U_{\b_1, \b_2,p}$, we obtain that $S_{\b_1,p}$ is bounded from $X(0,1)$ to $Y(0,1)$.

\emph{(iii) implies (iv).} Fix $f\in\M(0,1)$. First, note that
\begin{equation}\label{gaussible:calderon_thm:eq4}
T_{\b_1, \b_2, p} f(t) \lesssim \frac{f^*(\sigma(t^{1/p})^p)}{\b_1(\sigma(t^{1/p})^p)^p} + T_{\b_1, \b_2, p} f(\sigma^{-1}(t^{1/p})^p) \quad \text{for every $t\in(0,1)$}.
\end{equation}
Indeed, since $t \leq \sigma^{-1}(t^{1/p})^p$, we have that
\begin{align*}
    T_{\b_1, \b_2, p} f(t) &\leq \sup_{t\leq s\leq \sigma^{-1}(t^{1/p})^p} \frac{f^*(\sigma(s^{1/p})^p)}{\b_1(\sigma(s^{1/p})^p)^p} + \sup_{\sigma^{-1}(t^{1/p})^p \leq s < 1} \frac{f^*(\sigma(s^{1/p})^p)}{\b_1(\sigma(s^{1/p})^p)^p} \\ &\leq \frac{f^*(\sigma(t^{1/p})^p)}{\b_1(t)^p} + T_{\b_1, \b_2, p} f(\sigma^{-1}(t^{1/p})^p) \\ &\approx \frac{f^*(\sigma(t^{1/p})^p)}{\b_1(\sigma(t^{1/p})^p)^p} + T_{\b_1, \b_2, p} f(\sigma^{-1}(t^{1/p})^p).
\end{align*}

\noindent It follows from \eqref{gaussible:calderon_thm:eq4} that
\begin{align}
	\int_0^1 T_{\b_1, \b_2, p} f(t) g^*(t)\,dt &\lesssim \int_0^1 \left(\frac{f^*(\sigma(t^{1/p})^p)}{\b_1(\sigma(t^{1/p})^p)^p}\right) g^*(t)\,dt + \int_0^1 T_{\b_1, \b_2, p} f(\sigma^{-1}(t^{1/p})^p) g^*(t)\,dt \notag \\
		&= \mathcal I_1(g) + \mathcal I_2(g)\label{gaussible:calderon_thm:eq1}
\end{align}
for every $g\in\M(0,1)$. As for $\mathcal I_1(g)$, by the change of variables $t = \sigma^{-1}(\tilde{t}^{1/p})^p$, H\"older's inequality \eqref{ch1:ri:holder}, and $(iii)$, we have that
\begin{equation}\label{gaussible:calderon_thm:eq2}
\begin{split}
    \mathcal I_1(g) &\approx \int_0^1 f^*(t) (U_{\b_1, \b_2,p} (|g|^{1/p})(t))^p\,dt \leq \lVert f \rVert_{\big(Y^\frac 1p\big)'(0,1)} \lVert U_{\b_1, \b_2,p} (|g|^{1/p})\rVert_{Y(0,1)}^p \\ &\lesssim \lVert f \rVert_{\big(Y^\frac 1p\big)'(0,1)} \lVert |g|^{1/p}\rVert_{X(0,1)}^p  = \lVert f \rVert_{\big(Y^\frac 1p\big)'(0,1)} \lVert g\rVert_{X^\frac 1p(0,1)}.
\end{split}
\end{equation}

\noindent As for $\mathcal I_2(g)$, since the function $\frac1{\b_1}$ is equivalent to a quasiconcave function on $(0,1)$, it follows from \cite[Lemma~4.10]{EdMiMuPi:20} that

\begin{align*}
    \int_0^t T_{\b_1, \b_2, p} f(\sigma^{-1}(s^{1/p})^p) \,ds = \int_0^t \sup_{s \leq \tau < 1} \frac{f^*(\tau)}{\b_1(\tau)^p} \,ds \lesssim \int_0^t \left( \frac{f^*(\tau)}{\b_1(\tau)^p} \right)^*(s) \,ds
\end{align*}
for every $t\in(0,1)$. Hence, by virtue of Hardy's lemma \eqref{ch1:ri:hardy-lemma},
\begin{equation}\label{gaussible:calderon_thm:eq3}
    \mathcal I_2(g) \lesssim \int_0^1 \left( \frac{f^*(s)}{\b_1(s)^p} \right)^*(t) g^*(t)\,dt.
\end{equation}

Finally, by combining \eqref{gaussible:calderon_thm:eq1} with \eqref{gaussible:calderon_thm:eq2} and \eqref{gaussible:calderon_thm:eq3}, and using H\"older's inequality \eqref{ch1:ri:holder} and the boundedness of $S_{\b_1,p}$,  we obtain

\begin{align*}
    \|T_{\b_1, \b_2, p} f\|_{\big(X^\frac 1p\big)'(0,1)} &= \sup_{\lVert g \rVert_{X^\frac 1p(0,1)} \leq 1 } \int_0^1 T_{\b_1, \b_2, p} f(t) g^*(t)\,dt \lesssim \sup_{\lVert g \rVert_{X^\frac 1p(0,1)} \leq 1 } \left(\mathcal I_1(g) + \mathcal I_2(g)\right) \\
		&\lesssim \lVert f \rVert_{\big(Y^\frac 1p\big)'(0,1)} + \sup_{\lVert g \rVert_{X^\frac 1p(0,1)} \leq 1 } \int_0^1 \left( \frac{f^*(s)}{\b_1(s)^p} \right)^*(t) g^*(t)\,dt \\
		&= \lVert f \rVert_{\big(Y^\frac 1p\big)'(0,1)} + \sup_{\lVert g \rVert_{X^\frac 1p(0,1)} \leq 1 } \int_0^1  \frac{f^*(t)}{\b_1(t)^p}|g(t)|\,dt \\
		&\leq  \lVert f \rVert_{\big(Y^\frac 1p\big)'(0,1)} + \sup_{\lVert g \rVert_{X^\frac 1p(0,1)} \leq 1 } \int_0^1  \frac{f^{**}(t)}{\b_1(t)^p}|g(t)|\,dt \\
		&= \lVert f \rVert_{\big(Y^\frac 1p\big)'(0,1)} + \sup_{\lVert g \rVert_{X^\frac 1p(0,1)} \leq 1 } \int_0^1 f^*(t) S_{\b_1,p}(|g|^{1/p})(t)^p\,dt \\
		&\leq \lVert f \rVert_{\big(Y^\frac 1p\big)'(0,1)} + \lVert f \rVert_{\big(Y^\frac 1p\big)'(0,1)}\sup_{\lVert g \rVert_{X^\frac 1p(0,1)} \leq 1 }\|S_{\b_1,p}(|g|^{1/p})\|_{Y(0,1)}^p\\
		&\lesssim \lVert f \rVert_{\big(Y^\frac 1p\big)'(0,1)} + \lVert f \rVert_{\big(Y^\frac 1p\big)'(0,1)}\sup_{\lVert g \rVert_{X^\frac 1p(0,1)}\leq 1 }\||g|^{1/p}\|_{X(0,1)}^p  \\ & \approx \lVert f \rVert_{\big(Y^\frac 1p\big)'(0,1)}.
\end{align*}
Hence $T_{\b_1, \b_2, p}$ is bounded from $\big(Y^\frac 1p\big)'(0,1)$ to $\big(X^\frac 1p\big)'(0,1)$.

\

\emph{(iv) implies (i).} Since $T$ is $(p,\b_1,\b_2)$-gaussible, by virtue of Hardy's lemma \eqref{ch1:ri:hardy-lemma} we have that
\begin{align*}
    \int_0^1 (Tf)^*(s)^pg^*(s)\,ds &\leq \int_0^1 \left[(Tf)^*(s) \b_1(s)\right]^p \left(\sup_{s \leq \tau < 1} \frac{g^*(\tau)}{\b_1(\tau)^p}\right)\,ds \\
		&\lesssim \int_0^1 f^*(\sigma^{-1}(s^{1/p})^p)^p \b_1(s)^p\b_2(s)^{-p} \,  T_{\b_1, \b_2, p} g(\sigma^{-1}(s^{1/p})^p)\,ds \\
		&\approx \int_0^1 f^*(s)^p \, T_{\b_1, \b_2, p} g(s)\,ds
\end{align*}
for every $g\in\M(0,1)$. Hence, by using H\"older's inequality \eqref{ch1:ri:holder} on the right-hand side and $(iv)$,

\begin{equation*}
    \int_0^1 (Tf)^*(s)^pg^*(s)\,ds \lesssim \lVert |f|^p \rVert_{X^\frac 1p(0,1)} \lVert T_{\b_1, \b_2, p} g \rVert_{\big(X^\frac 1p\big)'(0,1)} \lesssim \lVert f \rVert_{X(R,\mu)}^p \lVert g \rVert_{\big(Y^\frac 1p\big)'(0,1)},
\end{equation*}
whence, by taking the supremum over all $g$ from the unit ball of $\big(Y^\frac 1p \big)'(0,1)$, we obtain that $T$ is bounded from $X(R,\mu)$ to $Y(S,\nu)$.

\end{proof}

It turns out that statement Theorem~\ref{gaussible:calderon_thm:restated}(iii) is often actually equivalent to (in turn, so are the other three statements):
\begin{itemize}
\item[(iii')]  The operator $U_{\b_1, \b_2,p}$ is bounded from $X(0,1)$ to $Y(0,1)$.
\end{itemize}
\begin{theorem}\label{gaussible:Usigma_bdd_implies_S_bdd}
Let $p\in(0,\infty)$ and $Y(0,1)$ be a $p$-convex r.i.~quasi-Banach function space. Let $(\b_1, \b_2) \in \mathcal B_p$. Furthermore, assume that the function $(0,1)\ni\tau\mapsto \b_1(\tau) \b_2(\tau)^{-1}$ is strictly decreasing and that
\begin{equation*}
 \lim_{s\to0^+}\b_2(s)^p\int_s^1 \frac{d\tau}{\tau \b_1(\tau)^p}\in(0,\infty).
\end{equation*}
We have that
\begin{equation}\label{gaussible:Usigma_bdd_implies_S_bdd:norm_ineq2}
\|S_{\b_1, p}(f^*)\|_{Y(0,1)}\lesssim \|U_{\b_1, \b_2,p} f\|_{Y(0,1)} \quad \text{for every $f\in\M(0,1)$}.
\end{equation}
Moreover, let $X(0,1)$ be another r.i.~quasi-Banach function space that is $p$-convex. If $U_{\b_1, \b_2,p}$ is bounded from $X(0,1)$ to $Y(0,1)$, so is $S_{\b_1,p}$.
\end{theorem}
\begin{proof}
First, note that $t< \sigma^{-1}(t^{1/p})^p$ for every $t\in(0,1)$ (recall Remark~\ref{gaussible:rem:properties_of_sigma}(ii)).

Next, since $Y^\frac 1p(0,1)$ is an r.i.~Banach function space, in order to prove \eqref{gaussible:Usigma_bdd_implies_S_bdd:norm_ineq2}, by virtue of the Hardy--Littlewood--P\'olya principle \eqref{ch1:ri:HLP} it is sufficient to show that
\begin{equation}\label{gaussible:Usigma_bdd_implies_S_bdd:eq1}
\int_0^t\int_s^1\frac{f^*(\tau)^p}{\tau\b_1(\tau)^p}\, d\tau\, ds\lesssim \int_0^t f^*(\sigma^{-1}(s^{1/p})^p)^p \b_2(s)^{-p} \, ds
\end{equation}
for every $t\in(0,1)$ and every $f\in\M(0,1)$ with a multiplicative constant independent of $f$ and $t$. Fix such $f$ and $t$. By Fubini's theorem, the left-hand side of \eqref{gaussible:Usigma_bdd_implies_S_bdd:eq1} is equal to
\begin{equation}\label{gaussible:Usigma_bdd_implies_S_bdd:eqi2}
\int_0^t\frac{f^*(s)^p}{\b_1(s)^p} \, ds + t\int_t^1\frac{f^*(s)^p}{s\b_1(s)^p} \, ds,
\end{equation}
and, by the change of variables $\tilde{s}=\sigma^{-1}(s^{1/p})^p$, the right-hand side of \eqref{gaussible:Usigma_bdd_implies_S_bdd:eq1} is equivalent to
\begin{equation}\label{gaussible:Usigma_bdd_implies_S_bdd:eqi3}
\int_0^{\sigma^{-1}(t^{1/p})^p}\frac{f^*(s)^p}{\b_1(s)^p}\, ds = \int_0^t \frac{f^*(s)^p}{\b_1(s)^p}\, ds + \int_t^{\sigma^{-1}(t^{1/p})^p} \frac{f^*(s)^p}{\b_1(s)^p}\, ds.
\end{equation}
In the light of \eqref{gaussible:Usigma_bdd_implies_S_bdd:eqi2} and \eqref{gaussible:Usigma_bdd_implies_S_bdd:eqi3}, in order to prove \eqref{gaussible:Usigma_bdd_implies_S_bdd:eq1}, it is sufficient to show that
\begin{equation}\label{gaussible:Usigma_bdd_implies_S_bdd:eq4}
t\int_t^1\frac{f^*(s)^p}{s\b_1(s)^p} \, ds \lesssim \int_t^{\sigma^{-1}(t^{1/p})^p} \frac{f^*(s)^p}{\b_1(\sigma(s^{1/p})^p)^p}\, ds
\end{equation}
with a multiplicative constant independent of $f$ and $t$. To this end, owing to the monotone convergence theorem and the fact that every nonnegative, nonincreasing function on $(0,1)$ is the pointwise limit of a nondecreasing sequence of nonincreasing simple functions on $(0,1)$, it is actually sufficient to prove \eqref{gaussible:Usigma_bdd_implies_S_bdd:eq4} for $f^*=\chi_{(0,a)}$, where $a\in(0,1)$. Therefore, \eqref{gaussible:Usigma_bdd_implies_S_bdd:eq4} will follow once we prove that
\begin{equation}\label{gaussible:Usigma_bdd_implies_S_bdd:eq5}
t\int_t^1\frac{\chi_{(0,a)}(s)}{s\b_1(s)^p} \, ds \lesssim \int_t^{\sigma^{-1}(t^{1/p})^p} \frac{\chi_{(0,a)}(s)}{\b_1(\sigma(s^{1/p})^p)^p}\, ds
\end{equation}
for every $a\in(0,1)$ with a multiplicative constant independent of $a$ and $t$. We claim that
\begin{equation}\label{gaussible:Usigma_bdd_implies_S_bdd:eq6}
s\int_s^1\frac{d\tau}{\tau\b_1(\tau)^p} \lesssim \int_s^{\sigma^{-1}(s^{1/p})^p} \frac{d\tau}{\b_1(\sigma(\tau^{1/p})^p)^p} \quad \text{for every $s\in(0,1)$}.
\end{equation}
Before we set out to prove the claim, we will make three observations. First, since the function $\tau\mapsto \b_1(\tau)^p \b_2(\tau)^{-p}$, $\tau\in(0,1)$, is decreasing, we have that
\begin{equation*}
\int_0^1\b_1(\tau)^p \b_2(\tau)^{-p}\d{\tau} > \lim_{s\to 1^-}\b_1(s)^p \b_2(s)^{-p}.
\end{equation*}
Second,
\begin{equation*}
\lim_{s\to0^+}\b_1(s)^p\int_s^1\frac{d\tau}{\tau\b_1(\tau)^p}=\infty \qquad \text{ and } \qquad \lim\limits_{s\to0^+}\frac{\b_1(s)}{\b_2(s)}=\infty,
\end{equation*}
for
\begin{equation*}
\frac{\b_1(s)^p}{\b_2(s)^p} \gtrsim \b_1(s)^p\int_s^1\frac{d\tau}{\tau\b_1(\tau)^p}\geq\int_s^1\frac{d\tau}{\tau} \quad \text{for every $s\in(0,1)$}.
\end{equation*}
Third, in order to prove \eqref{gaussible:Usigma_bdd_implies_S_bdd:eq6}, it is sufficient to prove that the inequality is valid near $0^+$ and near $1^-$ inasmuch as
\begin{equation*}
\sup_{s\in[c,d]}\frac{s\int_s^1\frac{d\tau}{\tau\b_1(\tau)^p} }{\int_s^{\sigma^{-1}(s^{1/p})^p} \frac{d\tau}{\b_1(\sigma(\tau^{1/p})^p)^p}}<\infty \quad \text{for every $0<c<d<1$}.
\end{equation*}
Set $M=(\int_0^1\b_1(\tau)^p \b_2(\tau)^{-p}\d{\tau})^{-1}$. As for the validity near $1^-$, note that, for every $s\in(0,1)$,
\begin{equation*}
\frac{s\int_s^1\frac{d\tau}{\tau\b_1(\tau)^p} }{\int_s^{\sigma^{-1}(s^{1/p})^p} \frac{d\tau}{\b_1(\sigma(\tau^{1/p})^p)^p}}\lesssim \frac{\frac{(1-s)}{\b_1(s)^p}}{\frac{\sigma^{-1}(s^{1/p})^p-s}{\b_1(s)^p}}=\frac{1-s}{\sigma^{-1}(s^{1/p})^p-s}
\end{equation*}
and that both numerator and denominator on the right-hand side goes to $0$ as $s\to1^-$. Hence, owing to L'H\^{o}pital's rule,
\begin{equation*}
\lim_{s\to1^-}\frac{1-s}{\sigma^{-1}(s^{1/p})^p-s}=\lim_{s\to1^-}\frac{-1}{M\b_1(s)^p \b_2(s)^{-p}-1}=\frac{-1}{M(\lim_{s\to 1^-}\b_1(s)^p \b_2(s)^{-p})-1} \in (0,\infty).
\end{equation*}
As for the validity near $0^+$, we use L'H\^{o}pital's rule again to obtain that
\begin{align*}
\lim_{s\to0^+}\frac{s\int_s^1\frac{d\tau}{\tau\b_1(\tau)^p} }{\int_s^{\sigma^{-1}(s^{1/p})^p} \frac{d\tau}{\b_1(\sigma(\tau^{1/p})^p)^p}} &= \lim_{s\to0^+}\frac{\int_s^1\frac{d\tau}{\tau\b_1(\tau)^p}  - \frac1{\b_1(s)^p}}{\frac{M \b_1(s)^p \b_2(s)^{-p}}{\b_1(s)^p} - \frac1{\b_1(\sigma(s^{1/p})^p)^p}}\\
&=\lim_{s\to0^+}\frac{\b_1(s)^p\int_s^1\frac{d\tau}{\tau\b_1(\tau)^p}  - 1}{M \b_1(s)^p \b_2(s)^{-p} - \frac{\b_1(s)^p}{\b_1(\sigma(s^{1/p})^p)^p}}\\
&=\frac1{M}\left( \lim_{s\to0^+}\b_2(s)^p\int_s^1 \frac{d\tau}{\tau \b_1(\tau)^p} \right) \in (0,\infty).
\end{align*}
Therefore, \eqref{gaussible:Usigma_bdd_implies_S_bdd:eq6} is valid. Having \eqref{gaussible:Usigma_bdd_implies_S_bdd:eq6} at our disposal, it is now easy to prove \eqref{gaussible:Usigma_bdd_implies_S_bdd:eq5}. If $a\leq t$, then \eqref{gaussible:Usigma_bdd_implies_S_bdd:eq5} plainly holds. If $t<a\leq \sigma^{-1}(t^{1/p})^p$, then
\begin{align*}
t\int_t^1\frac{\chi_{(0,a)}(s)}{s\b_1(s)^p} \, ds &= t\int_t^a\frac{ds}{s\b_1(s)^p} \leq \int_t^a\frac{ds}{\b_1(s)^p}\\
&\approx \int_t^a\frac{ds}{\b_1(\sigma(s^{1/p})^p)^p} = \int_t^{\sigma^{-1}(t^{1/p})^p}\frac{\chi_{(0,a)}(s)}{\b_1(\sigma(s^{1/p})^p)^p} \, ds.
\end{align*}
If $a>\sigma^{-1}(t^{1/p})^p$, then
\begin{align*}
t\int_t^1\frac{\chi_{(0,a)}(s)}{s\b_1(s)^p} \, ds &\leq t\int_t^1\frac{ds}{s\b_1(s)^p} \lesssim \int_t^{\sigma^{-1}(t^{1/p})^p} \frac{ds}{\b_1(\sigma(s^{1/p})^p)^p}\\
&= \int_t^{\sigma^{-1}(t^{1/p})^p}\frac{\chi_{(0,a)}(s)}{\b_1(\sigma(s^{1/p})^p)^p} \, ds,
\end{align*}
in which the multiplicative constant is that from \eqref{gaussible:Usigma_bdd_implies_S_bdd:eq6}. Hence \eqref{gaussible:Usigma_bdd_implies_S_bdd:eq5} is true, which completes the proof of \eqref{gaussible:Usigma_bdd_implies_S_bdd:norm_ineq2}.

Finally, let $X(0,1)$ be an r.i.~quasi-Banach function space that is $p$-convex and assume that $U_{\b_1, \b_2,p}\colon X(0,1)\to Y(0,1)$ is bounded. It follows from \cite[Corollary~9.8]{CiPiSl:15} (cf.~\cite[Theorem~1]{Pe:20}) that
\begin{align*}
\left\|\int_t^1\frac{|f(s)|^p}{s\b_1(s)^p}\, ds\right\|_{Y^\frac1{p}(0,1)} &\lesssim \||f|^p\|_{X^\frac1{p}(0,1)} \quad \text{for every $f\in\M(0,1)$}\\
\intertext{if and only if}
\left\|\int_t^1\frac{f^*(s)^p}{s\b_1(s)^p}\, ds\right\|_{Y^\frac1{p}(0,1)} &\lesssim \||f|^p\|_{X^\frac1{p}(0,1)} \quad \text{for every $f\in\M(0,1)$}.
\end{align*}
Owing to this equivalence, $S_{\b_1, p}\colon X(0,1)\to Y(0,1)$ is bounded if (and only if)
\begin{equation}\label{gaussible:Usigma_bdd_implies_S_bdd:eq7}
\|S_{\b_1, p}(f^*)\|_{Y(0,1)}\lesssim \|f\|_{X(0,1)} \quad \text{for every $f\in\M(0,1)$}.
\end{equation}
Thanks to \eqref{gaussible:Usigma_bdd_implies_S_bdd:norm_ineq2}, we have that
\begin{equation*}
\|S_{\b_1, p}(f^*)\|_{Y(0,1)} \lesssim \|U_{\b_1, \b_2, p}f\|_{Y(0,1)} \lesssim \|f\|_{X(0,1)} \quad \text{for every $f\in\M(0,1)$},
\end{equation*}
whence \eqref{gaussible:Usigma_bdd_implies_S_bdd:eq7} follows.
\end{proof}

We shall finish by illustrating our results with a particular example. Recall that the function $\ell\colon (0,1)\to(0,\infty)$ is defined as $\ell(t)=1-\log(t)$, $t\in(0,1)$. Set $\b_1 = \ell^\alpha$, $\b_2 = \ell^{-\beta}$. Let $p\in(0,\infty)$. It is a matter of straightforward computations to check that $(\b_1, \b_2)\in\mathcal B_p$ if and only if $\alpha, \beta\geq 0$ and either $\alpha + \beta\geq\frac1{p}$ and $\beta>0$ or $\alpha > \frac1{p}$ and $\beta=0$. Moreover, if either $0\leq\alpha<\frac1{p}$ and $\alpha + \beta = \frac1{p}$ or $\alpha>\frac1{p}$ and $\beta=0$, then the pair $(\b_1, \b_2)$ also satisfies the assumptions of Theorem~\ref{gaussible:Usigma_bdd_implies_S_bdd}. Therefore, by combining Theorems~\ref{gaussible:calderon_thm:restated} and \ref{gaussible:Usigma_bdd_implies_S_bdd}, we obtain the following important particular example. If $\alpha>\frac1{p}$, then $L^{\infty,\b_2}=L^\infty$, and so this case is not so interesting.

\begin{theorem}
\noindent Let $p\in(0, \infty)$ and $0\leq\alpha<\frac1{p}$. Set $\beta=\frac1{p}-\alpha$. Let $X(R,\mu)$ and $Y(S,\nu)$ be r.i.~quasi-Banach function spaces that are $p$-convex. The following four statements are equivalent.

\begin{enumerate}[(i)]
    \item Every $(p,\b_1,\b_2)$-gaussible operator $T$ is bounded from $X(R,\mu)$ to $Y(S,\nu)$.

    \item Every operator $T$ defined on $L^p(R,\mu)$ having values in $\M(S,\nu)$ that satisfies
		\begin{equation*}
			K(Tf, t; L^p(\log L)^\alpha, \exp L^\frac1{\beta}) \lesssim K(f,t;L^p,L^\infty) \quad \text{for every $f\in L^p(R,\mu)$ and $t\in (0,1)$}
		\end{equation*}
		is bounded from $X(R,\mu)$ to $Y(S,\nu)$.

    \item The operators $U_{\ell^\alpha, \ell^{-\beta},p}$ is bounded from $X(0,1)$ to $Y(0,1)$.
				
		\item The operator $T_{\ell^\alpha, \ell^{-\beta}, p}$ is bounded from $\big(Y^\frac 1p\big)'(0,1)$ to $\big(X^\frac 1p\big)'(0,1)$.

\end{enumerate}

\end{theorem}

\bibliography{Kinequality}

\end{document}